\documentclass{fic-l}

\newtheorem{theorem}{Theorem}[section]
\newtheorem{lemma}[theorem]{Lemma}
\newtheorem{corollary}[theorem]{Corollary}
\newtheorem{conjecture}[theorem]{Conjecture}
\newtheorem{proposition}[theorem]{Proposition}

\theoremstyle{definition}
\newtheorem{definition}[theorem]{Definition}
\newtheorem{example}[theorem]{Example}

\theoremstyle{remark}
\newtheorem{remark}[theorem]{Remark}

\numberwithin{equation}{section}

\usepackage{amscd}

\newcommand{\ra}{\rightarrow}
\newcommand{\sign}{\operatorname{sign}}
\newcommand{\ks}{\operatorname{ks}}
\newcommand{\kervaire}{\operatorname{kervaire}}
\newcommand{\C}{\mathbb{C}}
\newcommand{\R}{\mathbb{R}}
\newcommand{\Z}{\mathbb{Z}}

\begin{document}

\title{The Borel/Novikov conjectures and stable diffeomorphisms of 4-manifolds}

\author{James F. Davis}
\address{Department of Mathematics\\
Indiana University\\
Bloomington, IN 47401\\
USA}

\email{jfdavis@indiana.edu}

\thanks{Supported by the Alexander von Humboldt-Stiftung
and the National Science Foundation.  The author wishes to 
thank the Johannes Gutenberg-Universit\"at in Mainz for its hospitality while this
work was carried out.}

\date{December 1, 2004}

\begin{abstract}
Two 4-manifolds are stably diffeomorphic if they become diffeomorphic after connected sum with $S^2 \times S^2$'s. This paper shows that two closed, orientable, homotopy equivalent, smooth 4-manifolds are stably diffeomorphic provided a certain map from the second homology of the fundamental group with coefficients in $\Z_2$ to the L-theory of the group is injective. This injectivity is implied by the Borel/Novikov conjecture for torsion-free groups, which is known for many groups. There are also results concerning the homotopy invariance of the Kirby-Siebenmann invariant. The method of proof is to use Poincare duality in Spin bordism to translate between Wall's classical surgery and Kreck's modified surgery.
\end{abstract}

\maketitle

\section{Introduction}

Two smooth 4-dimensional manifolds $M$ and $N$ are {\em stably diffeomorphic} if for some
non-negative integers $r$ and $s$, the connected sum  
$M\# r(S^2\times S^2)$ is diffeomorphic to $N\# s(S^2\times S^2)$.  
This sort of stabilization plays a fundamental role in 4-dimensional
topology, see, for example \cite{Wall(64)}, \cite{Cappell-Shaneson}.  This paper examines the extent to which homotopy equivalent, smooth manifolds are stably diffeomorphic.  In this
paper we make the following two conjectures, relate them to standard 
conjectures in manifold theory, and thereby prove the following two conjectures for
large classes of fundamental groups.
 
\begin{conjecture}  \label{c1} If $M$ and $N$ are closed, orientable, smooth
4-manifolds which are homotopy equivalent and have torsion-free fundamental group, 
then they are stably diffeomorphic.
\end{conjecture}

\begin{conjecture}  \label{c2} If $M$ and $N$ are closed, orientable, topological 4-manifolds 
which are homotopy equivalent, have torsion-free fundamental group, and have the same
Kirby-Siebenmann invariant, then they are stably homeomorphic.
\end{conjecture}

In the simply-connected case the validity of Conjecture \ref{c1}  is well-known by the work of Wall 
\cite{Wall(64)}, who showed that homotopy equivalent, smooth, simply-connected 
4-manifolds are h-cobordant, and that h-cobordant, smooth, simply-connected manifolds are
stably diffeomorphic.  Using gauge theory, Donaldson \cite{Donaldson(87)} showed that
they need not be diffeomorphic.  In the simply-connected case, Conjecture \ref{c2} follows from 
the work of Freedman \cite{Freedman(82)}, with the stronger conclusion that the
manifolds are actually homeomorphic.  (The Kirby-Siebenmann invariant of a topological 4-manifold $M$ is a class $\ks (M) \in H^4(M;\Z_2)$
which vanishes if and only if $M \times \R$ admits a smooth structure.) 

The study of the stable diffeomorphism type of 4-manifolds divides into two cases, depending on whether the universal cover does or does not admit a Spin structure.  The case where the universal cover does not admit a Spin structure is much easier to analyze.  The following theorem (proved in Section 2) follows easily from results of Kreck. 

\begin{theorem}[Kreck]  Let $M$ and $N$ be closed, homotopy equivalent 4-manifolds whose universal covers do not admit a Spin structure.
\begin{enumerate}
\item Suppose $M$ and $N$ are smooth manifolds.  Then $M$ and $N$ are stably diffeomorphic if and only if  there is a homotopy equivalence $h : M \to N$ so that $w_1M = h^* w_1N$.
\item Suppose $M$ and $N$ are topological manifolds.  Then $M$ and $N$ are stably homeomorphic if and only if they have same Kirby-Siebenmann invariants and if there is a homotopy equivalence $h : M \to N$ so that $w_1M = h^* w_1N$.
\end{enumerate}
\end{theorem}

Note that, in particular, that any two homotopy equivalent, closed, orientable smooth 4-manifolds whose universal covers have a non-trivial second Stiefel-Whitney class are stably diffeomorphic. Note also that when the universal cover is not Spin, Conjectures \ref{c1} and \ref {c2} hold for any fundamental group. 

The case where the universal covers are Spin is much more subtle.  Teichner in his thesis \cite[Example 8.2.4]{Teichner(92)}
constructed an example of two closed, orientable, homotopy equivalent, smooth 4-manifolds 
with finite fundamental group which are not stably diffeomorphic.  

We put the question of whether two homotopy equivalent 4-manifolds are stably diffeomorphic in the context of surgery theory.  There is a map $\kappa_2 : 
H_2(\pi;\Z_2) \rightarrow L_4(\Z\pi)$, which appears in the surgery classification of high-dimensional manifolds.
(Here $L = L^h$, and refers to the Witt group of quadratic forms on free
$Z\pi$-modules.) 
As we shall see, this map is conjectured to be injective for all torsion-free
groups $\pi$, and has been shown to be injective in many cases, including
those listed in the corollary below.  In this note we give an alternate description of $\kappa_2$ and show:

\begin{theorem} \label{main}
If $\kappa_2 : H_2(\pi;{\Z}_2) \rightarrow L_4({\Z}\pi)$ is injective for a group $\pi$, 
and if $M$ and $N$ are closed, orientable, homotopy equivalent 4-manifolds with fundamental 
group isomorphic to $\pi$, then
\begin{enumerate}
\item If $M$ and $N$ are smooth manifolds, then they are stably diffeomorphic.
\item If $M$ and $N$ are topological manifolds with the same Kirby-Siebenmann invariant,
then they are stably homeomorphic.
\item If $M$ and $N$ are smooth manifolds with zero second Stiefel-Whitney
classes,
 then for some choice of Spin structure and 
identification of their fundamental groups with $\pi$, they are equal in $\Omega^{Spin}_4(B\pi)$.
\item If $M$ and $N$ are topological manifolds with zero second Stiefel-Whitney
classes,
 then for some choice of Spin structure and 
identification of their fundamental groups with $\pi$, they are equal in $\Omega^{TopSpin}_4(B\pi)$.
\item If $M$ and $N$ are topological manifolds so that the second Stiefel-Whitney classes of their
universal covers vanish, then $\ks (M) =  \ks (N)$.
\end{enumerate}
\end{theorem}

Our theorem holds even when the group is not torsion-free, 
but for finite groups  $\kappa_2$ need not be injective (see \cite[Prop. 7.4]{Hambleton-Milgram-Taylor-Williams}). 
All these results are well-known in the simply-connected case,
for example, part 5  
is clear, since the Kirby-Siebenmann invariant for Spin manifolds
is given by the signature divided by 8 considered modulo 2.  We conjecture that
the Kirby-Siebenmann invariant is a homotopy invariant for closed, orientable 4-manifolds with
torsion-free fundamental group whose universal cover is Spin, but this is false in general for manifolds
with finite fundamental group (see e.g. \cite[8.2.2]{Teichner(92)} or Example \ref{FKS} of this paper).  The techniques of our paper should also apply to non-orientable 4-manifolds; however we do not study them here.

 Theorem \ref{main} is the main result of this paper.  Its proof involves a comparison of C. T. C. Wall's surgery program
for studying homotopy equivalences with M. Kreck's surgery program for
studying stable diffeomorphisms.  We point out the ingredients of the proof.  Theorem \ref{Kreck} states Kreck's  classification of 4-manifolds up to stable diffeomorphism in terms of bordism.  Theorems \ref{SK} and \ref{SKno}  translate from surgery theory to bordism.  Corollary \ref{CSK} applies the characteristic class formulae \ref{charclass} to $\kappa_2$ to complete the proof.

We now switch to a discussion of the Borel/Novikov conjectures
and their relationships with $\kappa_2$.  Corollary \ref{cor} will apply the work of other mathematicians on the Borel/Novikov conjectures to Theorem \ref{main} to give a proof of the stable diffeomorphism conjectures in many cases.

A space is aspherical if its universal cover is contractible.  Let $M$ be a compact
aspherical  $n$-manifold.
\medskip

\noindent{\bf Borel Conjecture for $M$}{\bf :}  {\em Any homotopy equivalence from an $n$-dimensional manifold 
to $M$ which is a homeomorphism on the boundary is homotopic relative to the boundary
to a homeomorphism.}
\medskip

Equivalently the surgery structure set $\mathcal{S}(M\ \mbox{rel}\ \partial M)$ is trivial.  Reinterpreting
this in terms of the surgery exact sequence leads to:
\medskip

\noindent{\bf Strong Borel Conjecture for a torsion-free group $\pi$}{\bf :}  {\em The assembly
 map  \\ $H_*(\pi;\Bbb L _\cdot (\Bbb Z ))\rightarrow L_*(\Bbb Z \pi)$ is an isomorphism.}
\medskip

A more modest conjecture is that
this is an isomorphism when the Eilenberg-MacLane space  $B\pi$ is finite-dimensional.
Knowing the Borel conjecture for $M$, $M \times D^1$, $M \times D^2$, and $M \times D^3$, 
implies the strong Borel conjecture for  $\pi_1M$, as well as the vanishing of the Whitehead group  $\text{Wh}(\pi_1M)$.  Conversely, given a compact aspherical
manifold $M$ of dimension greater than four, the vanishing of the Whitehead group $\text{Wh}(\pi_1M)$ and  the strong Borel conjecture for $\pi_1M$
implies the Borel conjecture for $M \times D^i$, all $i$.  Farrell-Jones \cite{Farrell-Jones(93)} 
have proven the Borel conjecture
for closed manifolds with non-negative sectional curvature and their 
product with disks, provided the total dimension is greater than 4,
and hence the strong Borel conjecture is known for the fundamental groups of such manifolds.  
In addition the strong Borel conjecture can be proven for some groups 
built from the above groups using amalgamated products and HNN extensions, provided an extra condition
called square-root closed is satisfied (see Cappell \cite{Cappell}) and provided results
of Waldhausen \cite {Waldhausen} apply to show the 
Whitehead group is zero. 
The injectivity part of the strong Borel conjecture is often known as the integral Novikov
conjecture; the original Novikov conjecture is equivalent to that statement
that for all groups  $\pi$, the 
assembly map is a rational injection.  The integral Novikov conjecture is 
known for $\pi$ when $B\pi$ is a finite complex and $E\pi$ admits a certain type of
compactification by Carlsson-Pedersen \cite{Carlsson-Pedersen}.  (For example, this holds true when $\pi$ is a negatively-curved group in the sense
of Gromov or when $\pi$  is a discrete,
torsion-free, co-compact subgroup of a  virtually connected Lie group.)   There is a version of the above for non-orientable manifolds or equivalently for
groups equipped with an orientation character $w : \pi \ra \{\pm 1\}$, but the reader is advised
to consider only the orientable case $w = 1$ on a first reading.  The study of Borel/Novikov 
conjectures is a rapidly advancing area of mathematics.  For background the reader can consult
reports to three international congresses \cite{Farrell-Jones},
\cite{Quinn}, \cite{Hsiang} and references given there.  

The $i$-th homotopy group of the simply-connected $L$-spectrum $\Bbb L_\cdot (\Bbb Z)$ is
zero for $i$ odd, $\Bbb Z$ for $i = 4k$, and $\Bbb Z_2$ for $i = 4k+2$.
There is a natural splitting of $\Bbb L_\cdot (\Bbb Z)_{(2)}$ as a
wedge of Eilenberg-MacLane spectra {\cite{Taylor-Williams}}, and this gives an injection  
$H_2(B\pi;{\Bbb Z}_2) \rightarrow
H_4(B\pi;\Bbb L _\cdot (\Bbb Z ))$ and the composite with the assembly map is
called $\kappa_2$.  Thus $\kappa_2$ is conjecturally injective for all torsion-free groups and is
 injective whenever the integral Novikov conjecture holds.

\begin{corollary} \label{cor} If $\pi$ is
\begin{enumerate}
\item the fundamental group of a closed Riemannian manifold with sectional
curvature everywhere less than or equal to zero, or
\item a discrete, torsion-free subgroup of a virtually connected Lie group $G$ so that
$\pi\backslash G/K$ is compact where $K$ is a maximal compact subgroup of $G$, or
\item a negatively curved group in the sense of Gromov,
\end{enumerate}
then $\kappa_2$ is injective, and hence
for closed, orientable, homotopy equivalent 4-manifolds with
fundamental group isomorphic to $\pi$, the stable diffeomorphism conjectures
1-5 of Theorem \ref{main} hold true.
\end{corollary}

It would be interesting to give a purely algebraic definition of $\kappa_2$.  It would be more interesting to give four-dimensional proofs of the stable diffeomorphism 
conjectures, even in known cases, as this would shed light on the Borel/Novikov conjectures.

I would like to thank my colleagues at Mainz: Fang Fuquan, Matthias Kreck, and Wolfgang L\"uck, the
former for many discussions and the latter two for suggesting that I investigate the relationship
between the Borel/Novikov conjectures and four-dimensional topology.  I would also like to thank Qayum Khan for pointing out a subtlety in the proof of Proposition \ref{charclass}.

\section{Stable diffeomorphism of 4-manifolds}

M. Kreck \cite{Kreck} showed that stable diffeomorphism is a bordism question, and P. Teichner
\cite{Teichner(92)}, \cite{Teichner(93)} studied this question in dimension 4.  We review this
theory and use it to motivate our key invariant for homotopy equivalent 4-manifolds,
the dimension 2 Spin-Kervaire invariant.

Let $\xi:B\ra BO$  be a fibration\footnote{Throughout this paper there are parallel theories
for smooth manifolds and topological manifolds, with $BO$ replaced by $BTOP$.
In our theorems we will state both cases, but in the definitions and proofs
we will only deal with the smooth case unless the topological case involves a substantive
difference.}.  Let  $\Omega_*(\xi)$ be bordism classes of smooth manifolds 
equipped with a lift of the stable normal bundle $\nu$ to $\xi$.  (For details see \cite{Stong}
or \cite{Switzer}.)
Elements of $\Omega_*(\xi)$ are represented by pairs  $(M,\tilde{\nu}:M\ra B)$, where
$\xi\circ\tilde{\nu}=\nu:M \ra BO$.

If $\tilde{\nu}: M \ra B$ is a $k$-equivalence (i.e. $\pi_iM \ra \pi_iB$ is an isomorphism
for $i<k$ and a surjection for $i=k$), then $(M,\tilde{\nu})$ is called a {\em normal 
$(k-1)$-smoothing in $\xi$}.  Two normal $(k-1)$-smoothings $(M,\tilde{\nu})$ and
$(M^\prime,\tilde{\nu}^\prime)$ in $\xi$ are {\em diffeomorphic} if there is a diffeomorphism
$f : M \ra M^\prime$ so that $\tilde{\nu}$ and $\tilde{\nu}^\prime \circ f$ are  homotopic.
Let $NSt_{2n}(\xi)$ be stable diffeomorphism classes of $2n$-dimensional normal $(n-1)$-smoothings
in $\xi$.

The following theorem is a consequence of Section 4 of \cite{Kreck}.

\begin{theorem}[Kreck] \label{Kreck} Let $\xi:B \ra BO$ be a fibration  where the 
$n$-skeleton of $B$ has the homotopy type of a connected finite complex.  For $n\ge 2$, the following map
is a bijection  $$NSt_{2n}(\xi)\ra\Omega_{2n}(\xi).$$  For a fibration  $\xi:B \ra BTOP$
with $B$ as above and $n\ge 2$, the following is a bijection: 
$$NSt^{TOP}_{2n}(\xi)\ra\Omega^{TOP}_{2n}(\xi).$$
\end{theorem}

It is possible for the same manifold to represent two different elements of these sets by
composing one representation with an automorphism of $\xi$.  Let Aut$(\xi)$ be the group of
fiber homotopy equivalences of $\xi$; elements are given by homotopy equivalences
$f:B\ra B$ where $\xi\circ f = \xi$.  Aut$(\xi)$ acts on NSt$_{2n}(\xi)$ and 
$\Omega_{2n}(\xi)$ by $f(M,\tilde{\nu}) = (M,f\circ\tilde{\nu})$.

\begin{definition} {\em If $\xi:B\ra BO$ is $k$-coconnected (i.e. i.e. $\pi_iB \ra \pi_iBO$ is an isomorphism
for $i>k$ and a injection for $i=k$), then $\xi$ is {\em $(k-1)$-universal}. If $\tilde{\nu}:M \ra B$ is a normal $(k-1)$-smoothing in $\xi$ where 
$\xi$ is $(k-1)$-universal, then  $\xi : B \ra BO$ is the {\em normal $(k-1)$-type of $M$}. }
\end{definition}

Obstruction theory and the Moore-Postnikov factorization show that the 
normal $(k-1)$-type of $M$ exists and is unique up to fiber homotopy
equivalence.  Thus if $\xi$ is $(n-1)$-universal, $n \ge 2$ there is a bijection between stably diffeomorphism classes of manifolds with
normal $(n-1)$-type $\xi$ and $\Omega_{2n}(\xi)/$Aut$(\xi)$.

For stable diffeomorphism classes of 4-manifolds there are basically
two cases:  \linebreak $w_2(\widetilde{M}) \not= 0$ and $w_2(\widetilde{M}) = 0$.  The first case is quite simple; the following lemma implies that in this
case, homotopy equivalence implies stable diffeomorphism, even for 
manifolds with torsion in the fundamental group.

We assume all manifolds in this paper are connected.
If $M$ is a closed manifold $n$-manifold, then a {\em fundamental class} $[M] \in H_n(M;{\Bbb Z}^w)$
is a choice of generator for this infinite cyclic group, where ${\Bbb Z}^w$ denotes the 
${\Bbb Z}\pi_1M$-module given by the integers twisted by the orientation character $w = w_1M : \pi_1M \ra 
\{ \pm 1 \}$.  The following lemma is a easy consequence of Theorem \ref{Kreck} (cf. \cite{Teichner(92)}).

\begin{theorem}  \label{nospin}Two closed, smooth, orientable (resp. non-orientable)
4-manifolds $M$ and $N$ whose universal covers are not Spin are stably diffeomorphic iff there are fundamental classes $[M]$ and $[N]$ and an isomorphism $f:\pi_1M \to \pi_1N $ so  that    $w_1 M = w_1N  \circ f$,  $f_*[M] = [N]$, and  $\sign M = \sign N$ (resp. $\chi(M) \equiv \chi(N)$ mod 2).   Two
closed, topological, orientable (resp. non-orientable) 4-manifolds $M$ and $N$ 
whose universal covers are not Spin  are stably homeomorphic iff there are fundamental classes $[M]$ and $[N]$ and an isomorphism $f:\pi_1M \to \pi_1N $ so  that    $w_1 M = w_1N  \circ f$,  $f_*[M] = [N]$,   $\sign M = \sign N$ (resp. $\chi(M) \equiv \chi(N)$ mod 2), and 
$\ks (M) = \ks (N)$.
\end{theorem}

\begin{proof}
   Of the four cases, we will only write the proof for the 
simplest (smooth and orientable) and the hardest (topological and non-orientable).
In the smooth and orientable case 
the normal 1-type of $M$ is given by $\xi:BSO \times B\pi
\ra BO$, where $\xi$ is projection on the first factor followed by the 
double cover.  Thus  $\Omega_*(\xi) = \Omega^{SO}_*(B\pi)$.  The
A-H-S-S (= Atiyah-Hirzebruch Spectral Sequence) gives
$H_p(B\pi;\Omega^{SO}_q) \Rightarrow \Omega^{SO}_{p+q}(B\pi)$.  Since 
$\Omega^{SO}_q = {\Bbb Z},0,0,0,{\Bbb Z}$ for $q = 0,1,2,3,4$, we see
$\Omega^{SO}_4(B\pi) = \Omega^{SO}_4 \oplus H_4(\pi)$ and the result follows.

In the topological, non-orientable case, the normal 1-type of $M$ is given by \linebreak $\xi:BSTOP \times B\pi
\ra BTOP$, where $\xi$ classifies the product of the universal bundle $\gamma$ over $BSTOP$ with the orientation
line bundle $\zeta$ over $B\pi$.  We write $\Omega^{STOP}_*(B\pi;w) = \Omega_*(\xi)$. 
This is isomorphic to $\Omega^{STOP}_{*+1}(D(\zeta),S(\zeta))$, where the isomorphism is
given by pulling back the line bundle and the inverse map by taking the
transverse inverse image  of the zero section under 
a representative map $W \ra D(\zeta)$.   There is an A-H-S-S
 $H_p(B\pi;(\Omega^{STOP}_q)^w) \Rightarrow \Omega^{STOP}_{p+q}(B\pi;w)$.  Now
$\Omega^{STOP}_q = {\Bbb Z},0,0,0,{\Bbb Z} \oplus {\Bbb Z}_2$ for $q = 0,1,2,3,4$, 
with the extra class given by the Kirby-Siebenmann invariant \cite[p. 322]{Kirby-Siebenmann}.
 One deduces 
$$\Omega^{STOP}_4(B\pi;w) = {\Bbb Z}_2 \oplus {\Bbb Z}_2  \oplus H_4(B\pi;{\Bbb Z}^w)$$
with the invariants given by $\chi$ mod 2, ks, and the fundamental class. The result follows.
\end{proof}

The more interesting case is when the universal cover is Spin.  If M itself
is Spin, then the normal 1-type of $M$ is $\xi:BSpin \times B\pi \ra
BO$ and hence $\Omega_*(\xi) = \Omega^{Spin}_*(B\pi)$.  Thus two
Spin 4-manifolds with fundamental group isomorphic to $\pi$ are
stably diffeomorphic if for some choice of orientations, Spin structures,
and identification of the fundamental group, they represent the same
element in $\Omega^{Spin}_4(B\pi)$.  There is an A-H-S-S
$H_p(B\pi;\Omega^{Spin}_q) \Rightarrow \Omega^{Spin}_{p+q}(B\pi)$ and
$\Omega^{Spin}_q = {\Bbb Z},{\Bbb Z}_2,{\Bbb Z}_2,0,{\Bbb Z}$ for $q = 0,1,2,3,4$.

Teichner \cite{Teichner(93)} constructs the ``James spectral sequence''
to say that much the same is true when one only has $w_2(\widetilde{M})=0$.
The sequence
$$0 \ra H^2(\pi;{\Bbb Z}_2) \ra H^2(M;{\Bbb Z}_2) \ra H^2(\widetilde{M};{\Bbb Z}_2)$$ is exact and 
we let $w_i \in H^i(\pi;{\Bbb Z}_2), i = 1,2$ denotes the elements which
maps to $w_i(M)$.  The normal 1-type of $M$ is given as the
homotopy pullback

$$\begin{CD}  B @>>> B\pi \\
@VVV @VV w_1 \times w_2 V \\
BO @> w_1(\gamma) \times w_2(\gamma) >> K(\Z_2,1) \times K(\Z_2,2) \\
\end{CD}$$

(Homotopy pullback means convert the bottom and right hand maps to fibrations
and take the ordinary pullback.)  Denote the left hand map by
$\xi(B\pi,w_1,w_2)$ and the upper left hand space by $B(B\pi,w_1,w_2)$.  Use analogous notation when $B\pi$ is replaced by an arbitrary space $X$.

\begin{theorem}[\cite{Teichner(93)}]  Given classes $w_i \in H^i(X;{\Bbb Z}_2), i =1,2$, there is a spectral sequence  $H_p(X;(\Omega^{Spin}_q)^{w_1}) \Rightarrow 
\Omega_{p+q}(\xi(X,w_1,w_2))$.  
There is an analogous spectral sequence for topological manifolds.
\end{theorem}

Since $\Omega^{Spin}_3 = 0$,  $E^2_{2,2} = H_2(X;{\Bbb Z}_2) \ra E^\infty _{2,2}$ is surjective.
In the next section we will show that if $f : M \ra N$ is a homotopy
equivalence between closed 4-manifolds whose universal covers
have zero second Stiefel-Whitney classes, one can choose lifts $\tilde{\nu}_M$,
$\tilde{\nu}_N$ (``Spin structures'') so that $(M,\tilde{\nu}_M) -
(N,\tilde{\nu}_N)$ lies in $F_{2,2}$ and that the corresponding element
in $H_2(B\pi;{\Bbb Z}_2)/(d_2,d_3)$ is independent of the choice of lifts.
We call this element the {\em dimension 2 Spin-Kervaire invariant} of 
the homotopy equivalence.  Although we could give a direct proof that
this is well-defined, we prefer to identify it with the {\em codimension
2 Kervaire invariant} coming from surgery theory.  We do so in the next section.

There is a special case where the James spectral sequence can be replaced by a more familiar spectral sequence (cf. \cite[p. 53]{Teichner(93)}).

\begin{lemma}  \label{James} Let $\eta : X \to BO(k)$.
\begin{enumerate}

\item  The bundle $\xi(X, w_1(\eta),w_2(\eta))$ can be identified with 
$$
\xi :  BSpin \times X \to BO,
$$
classifying the product of  the universal Spin bundle over $BSpin$ and the bundle $\eta$.

\item  $\Omega_n(\xi)$ is given by bordism of closed smooth manifolds  $(M^n, f: M^n \to X)$ together with a Spin structure on $\nu_M \oplus f^*\eta$.

\item  $\Omega^{Spin}_{*+k}(D(\eta),S(\eta)) \cong \Omega_*(\xi)$.  The map from left to right is given by taking the transverse inverse image of the zero section.  The map from right to left is given by taking the pullback  $f^*(D(\eta),S(\eta))$.

\item  There is an isomorphism of the A-H-S-S of the pair with the James spectral sequence of $\Omega_*(\xi)$ from $E^2$ on.  The map on $E^2$ is the Thom isomorphism.

\end{enumerate}
\end{lemma}

\section{Surgery}
\begin{definition} {\em \cite{Wall(70)}, \cite{Browder} A {\em degree one normal map} is a map  $f :M \ra N$ between
closed manifolds equipped with fundamental classes
so that $f_*[M] = [N]$, together with a bundle
$\xi$ over $N$ and a stable trivialization of $\tau_M \oplus f^*\xi$.}
\end{definition}
\begin{example} {\em If $M$ is a closed, oriented $n$-manifold with a framing of the
stable normal bundle, let $f : M \ra S^n$ be a degree one map and
$\xi$ the trivial bundle.}
\end{example}  
\begin{example}
{\em Let $f:M \ra N$ be homotopy equivalence of closed manifolds,
together with a fundamental class for $M$.
 Choose a homotopy inverse $g$ of $f$ and let $\xi$ be $g^*\nu_M$.
Then a homotopy from $g\circ f$ to the identity gives an isomorphism of
$f^*\xi$ with $\nu_M$, and hence a framing of $\tau_M \oplus f^*\xi$.}
\end{example}
\begin{definition} {\em (See \cite[Theorem 2.23]{Milgram-Madsen}, also \cite[Section II.4]{Browder} for a more homotopy
theoretic approach.)
If $f$ is a degree one normal map, let $\lambda$ be the stable vector bundle $\nu_N - \xi$, 
and perturb the map $f$ from $M$ to the zero section to an embedding $M \ra E(\lambda)$.  
The stable trivialization leads to an extension to a codimension zero embedding 
$M \times D^k \ra E(\lambda)$
and by the Pontryagin-Thom construction to a map $T(\lambda) \ra \Sigma^kM_+$,
where $\Sigma^kM_+$ denote the $k$-fold reduced suspension of the disjoint union of
$M$ and a base point.  There is thus
a fiber homotopy equivalence  $S(\lambda \oplus 1) \ra N \times S^k$  where the map to the
first factor is given by projection to $N$ and the map to the second by mapping to the Thom space,
 composing with the Pontryagin-Thom map, followed by collapsing $M$ to a point.  A fiber
homotopy equivalence is classified by a map to $G/O =$ homotopy fiber$(BO \ra BG)$ where $G$ is the
monoid of self-homotopy equivalences of the sphere and $BG$ classifies spherical fibrations.
The homotopy class of the classifying map in $[N,G/O]$ is the {\em normal invariant of f}.}
\end{definition}

If two degree one normal maps to $N$ have the same normal invariant, then the maps are
{\em normally bordant}.  Geometrically, this means that the two bundles over $N$ are 
stably isomorphic
and that one can do surgery to get from one map to the other.  In the above two examples, there
is a well-defined normal bordism class, independent of choices.  Note that  $G/O$ and $G/TOP$
are connected, simply-connected and have $\pi_2 \cong {\Bbb Z}_2$.  This uses  the fact that
the homotopy groups of $G$ are the stable homotopy groups of spheres \cite[Corollary 3.8]{Milgram-Madsen}.  For a closed 2-manifold $N$,
not necessarily orientable, the bijection $[N,G/TOP] \ra {\Bbb Z}_2$ is called the Kervaire invariant.
A homotopy theoretic definition is given in \cite[Section III.4]{Browder}, a geometric one in \cite{Wall(70)}.

\begin{definition}
{\em If $f : M \ra N$ is a degree one normal map, the {\em codimension 2 Kervaire invariant} $\kervaire^2(f)\in
H^2(N;{\Bbb Z}_2)$ is the first obstruction to the normal invariant being null-homotopic.
Geometrically this corresponds to the cohomology class given by representing a homology class
by a 2-dimensional submanifold $P$ of $N$, making $f$ transverse to $P$ and assigning to 
$[P] \in H_2(N; \Z_2)$ the classical Kervaire invariant of the 2-dimensional normal map $f^{-1}(P) \ra P$.}
\end{definition}

The fundamental theorem of 
surgery theory \cite{Wall(70)} states that there are 
4-periodic abelian groups $L_n({\Bbb Z}\pi;w)$, 
natural in $(\pi;w)$, and a function $\theta : [N,G/O] \ra L_n({\Bbb Z}\pi;w)$ so that for $n > 4$,
($n$ = dimension $N$, $w = w_1(N)$), $\theta(\hat{f}) = 0$ if and only if $f$ is normally bordant to a homotopy
equivalence.  
The analogous thing is true in the topological category by \cite{Kirby-Siebenmann}.  

\begin{proposition}\label{charclass} Let $\pi$ be a group with orientation character $w$.  There are homomorphisms
$$I_0:H_0(\pi;{\Bbb Z}^w) \ra L_4({\Bbb Z}\pi;w)$$
$$\kappa_2:H_2(\pi;{\Bbb Z}_2)\ra L_4({\Bbb Z}\pi;w)$$  
so that for any degree one normal map $f : M \ra N$ between closed, connected 4-manifolds with classifying map
$\hat{f}: N\ra G/TOP$, and map $F : N \ra B\pi$ so that $w \circ F_*$ is the orientation character
of $N$, then one has the following characteristic class formulae:
\begin{align}
I_0((\sign M - \sign N)/8) + \kappa_2(F_*(\kervaire^2(f) \cap [N])) = 
\theta(\hat f)  \tag{$w=1$}  \\
I_0(\ks(M) - \ks(N) +\kervaire^2(f)^2) + \kappa_2(F_*(\kervaire^2(f) \cap [N])) = 
\theta(\hat f)  \tag{$w\not=1$} 
\end{align}
where we have identified $H_0(\pi;{\Bbb Z})$ with $\Bbb Z$ and when $ w\not = 1$ we have
identified $H^4(N;{\Bbb Z}_2)$ with $H_0(\pi;{\Bbb Z}^w)$.
\end{proposition}

\begin{proof}  By allowing surgeries in both the domain and range, Sullivan-Wall 
\cite[13B.3]{Wall(70)} show that the surgery obstruction map factors through a homomorphism
$$\theta:{\Omega}_*^{STOP}(B\pi \times G/TOP, B\pi \times pt;w) \ra L_*({\Bbb Z}\pi;w).$$ 
For $*=4$, the A-H-S-S identifies the domain of this map with 
$$
H_4(B\pi \times G/TOP, B\pi \times pt;{\Bbb Z}^w),
$$ 
with the identification given by the 
image of the fundamental class.
  The fourth stage of the Postnikov tower for $G/TOP$
is $K({\Bbb Z}_2,2) \times K(\Bbb Z,4)$ (see \cite[p. 329]{Kirby-Siebenmann}).  Recall
$H_4(K(\Z_2,2);\Z_2) = \Z_2$. We thus have a 
splitting of the domain of $\theta$ into 
$$H_0(\pi;{\Bbb Z}^w) \oplus H_0(\pi; \Z_2) \oplus
H_2(\pi; {\Bbb Z}_2)$$
 and this defines maps  $I_0$, $J_0$, and  $\kappa_2$ so that $\theta = I_0 \oplus J_0 \oplus \kappa_2$.  

The generators of the low-dimensional cohomology of $G/TOP$ are given by $k \in H^2(G/TOP;{\Bbb Z}_2)$
and $l \in H^4(G/TOP;\Bbb Z)$.  For a degree one normal map $f : M \ra N$ between {\em oriented}
4-manifolds with classifying map $\hat{f}$, one has 
\begin{align}
\hat{f}^*(k) &= \kervaire^2(f)  \notag \\
\hat{f}^*(l) &=((\sign M - \sign N)/8) \mu_N \notag
\end{align}
where $\mu_N \in H^4N$ is dual to the fundamental class. It follows that
$$
\theta(\hat f) = I_0((\sign M - \sign N)/8) + J_0(\kervaire^2(f)^2)  +\kappa_2(F_*(\kervaire^2(f) \cap [N])) 
$$
where we have identified $H_0(\pi)$ with $\Z$ and have identified $H_0(\pi;\Z_2) = \Z_2 = H^4(N;\Z_2)$.

There is a topological manifold $Ch^4$ homotopy equivalent, but not homeomorphic to $\C P^2$ (see \cite {Freedman(82)}).  The corresponding homotopy equivalence $h$ has non-trivial normal invariant, hence non-trivial codimension 2 Kervaire invariant.  The above formula then shows that $J_0(\kervaire^2(h)^2) = 0$.  Hence, by naturality in $(\pi;w)$, the homomorphism $J_0$ is always zero.

If $f: M \ra N$ is a degree one normal map between 4-manifolds with $N$ {\em non-orientable}, then 
by \cite[p. 329]{Kirby-Siebenmann}, one has
\begin{align}
\hat{f}^*(k) &= \kervaire^2(f)  \notag \\
\hat{f}^*(l) &= \kervaire^2(f)^2 + (\ks M - \ks N) \notag
\end{align}
after identifying $H^4(N;\Bbb Z)$ with $\Bbb Z_2$.  The characteristic class formula
follows.

\end{proof}

\begin{remark}
 This proposition is a special case of the characteristic class formulae of
\cite{Taylor-Williams}, although the proof is considerably simpler in dimension 4.  The 
$|E_8|$-manifold shows that in the simply-connected case that $I_0$ is an isomorphism.
In general $I_0(1)$ is the image of the $E_8$-quadratic form.  Since $1 \ra \pi$, is
a split injection of groups, in the oriented case the map $I_0$ is always injective.
In the non-oriented case this need not be the case.  (See \cite{Teichner(97)} for geometric applications.)
\end{remark}

We wish to identify the codimension 2 Kervaire invariant with the dimension 2 Spin-Kervaire invariant.
In doing so we will show that the dimension 2 Spin-Kervaire invariant is well-defined and 
draw the connection between stable diffeomorphism and $\kappa_2$ and thereby prove our
theorem.  The problem is that one invariant lies in a cohomology theory and 
 and the other lies in a homology theory, so we are led naturally to duality.
Recall that for any generalized homology theory there is a generalized cohomology theory
\cite{Whitehead} and that if a manifold is orientable with respect to that theory there
is a Poincar\'e duality between them.  For bordism and cobordism there is a 
geometric interpretation \cite{Atiyah}, \cite{Conner-Floyd}.  For a compact, Spin manifold
$N$ of dimension $n$, there is a duality  $\Omega^{Spin}_{k}(N,\partial N) \cong
\Omega_{Spin}^{n-k}(N)$.
Geometrically this goes as follows:  if  $k < n/2$, represent an element of Spin bordism by
an embedding $(V,\partial V) \ra (N, \partial N)$.
  Since $V$ and $N$ both have Spin structures there is a Spin 
structure on the normal bundle of $V$ in $N$.  
The Pontryagin-Thom construction gives a map $N \ra T(\nu_{V \subset N}) \ra T(\gamma_{n-k})$,
and hence a map  $N \ra MSpin(n-k)$.
If $k$ is large, we suspend by embedding $V$ in  $N \times {\Bbb R}^p$, and the Pontryagin-Thom
construction gives a map $\Sigma^pN_+ \ra MSpin(n-k+p)$, which gives a representative of 
Spin cobordism.
The inverse map is given by restricting a map $\Sigma^pN_+ \ra MSpin(n-k+p)$ to $N \times {\Bbb R}^p$, making this transverse to the zero section $BSpin(n-k+p)$, and 
letting $V$ be the inverse image.

\begin{lemma} \label{duality lemma} If $N$ is a compact Spin manifold of dimension $n$, 
then the A-H-S-S's for $\Omega_{Spin}^*(N)$ and $\Omega^{Spin}_*(N,\partial N)$ are dual in the sense
that $\cap [N] : (E^{p,q}_r,d_r) \ra (E^r_{n-p,-q},d^r)$ is an isomorphism for all $r \ge 2$,
and the induced map on $E_\infty$ corresponds to Spin bordism Poincar\'e duality.
\end{lemma}
\begin{proof} Poincar\'e duality can be factored as the composite of a Thom isomorphism map and 
a Spanier-Whitehead duality isomorphism.  (See \cite{Switzer} and \cite{Browder}.)
The Thom isomorphism goes from the cohomology of $N$ to the reduced cohomology of the Thom
space of the normal bundle $\nu$ of $N$ and S-duality maps further to the reduced homology of
$N/\partial N$.  This works equally well for Spin bordism and ordinary homology.  The Thom
classes and fundamental classes correspond under the natural transformation  MSpin $\ra H(
{\Bbb Z})$.

By choosing the cell structure on $T(\nu)$ corresponding to that of $N$, one can even guarantee
that $\Omega^k_{Spin}(N) \ra \widetilde{\Omega}^{n+k}_{Spin}(T(\nu))$ induces an isomorphism on
the A-H-S-S's from $E_1$ on, which on $E_2$ is given by cupping with the Thom class.  Next suppose
$X$ is a finite complex and we have a S-duality map $\mu: X \wedge X^* \ra S^L$.  By restriction
we have a map $X^i \wedge (X^*)^{L-i} \ra S^L$ and hence a map $\Omega^k_{Spin}(X^i) \ra
\Omega^{Spin}_{L-k}((X^*)^{L-i})$.  Since S-duality is well-behaved with respect to cofibrations
(see, for example, \cite[14.31]{Switzer}), it follows that there are maps from the Spin cobordism exact sequence
of the pair $(X^i,X^{i-1})$ to the Spin bordism exact sequence of the pair $((X^*)^{L-i},
(X^*)^{L-i-1})$.  This gives a map of A-H-S-S's for $\Omega^k_{Spin}(X)$ and $\Omega^{Spin}_{L-k}(X^*)$
from $E_1$ on, which is the S-duality isomorphism both at the $E_2$ level and
the Spin bordism level.  Since $T(\nu) (= X)$
and $N/\partial N$ are S-dual, this completes the proof.    
\end{proof}

We omit the proof of the following lemma.

\begin{lemma} \label{obstruction lemma} 
Let $a:\Sigma^{-p-q}X_+ \ra MSpin$ be a stable, base point preserving map which is trivial on the
$(-q-1)$-skeleton.  Then the image of $a$ in $E^{p,q}_\infty$ equals the image of
the first obstruction to null-homotopy $\vartheta(a) 
\in H^p(X;\Omega^q_{Spin})$.  Here we consider $E^{p,q}_\infty$ as both a subquotient of
$\Omega^{p+q}_{Spin}(X)$ and of $H^p(X;\Omega^q_{Spin})$.
\end{lemma}

If two maps are normally bordant, then their domains are stably diffeomorphic by \cite{Kreck-Schafer}
or \cite[Section 4]{Kreck}.  This leads one to expect that the obstruction to two maps being normally
bordant should determine whether the domains are stably diffeomorphic.  The precise relationship is
given by the
following theorem which identifies the dimension 2 Spin-Kervaire invariant with the image of the 
codimension 2 Kervaire invariant.  Technically speaking, it is the main result of this paper.
For simplicity we first consider the Spin case, and later
briefly indicate the changes necessary for non-Spin manifolds.

\begin{theorem} \label{SK} Let $f:M\ra N$  be a degree one normal map between closed, smooth 4-manifolds 
which both have the same signature and both
admit Spin structures.  Given any Spin structure $\tilde{\nu}_M :M \ra BSpin$ there exists
a unique Spin structure $\tilde{\nu}_N : N \ra BSpin$ so that $\alpha = [(M,  \tilde{\nu}_M),f] -
[(N,  \tilde{\nu}_N), Id ]$ is in the filtration subgroup  $F_{2,2}
(\Omega^{Spin}_4(N)) = Im(\Omega^{Spin}_4(N^{(2)}) \ra \Omega^{Spin}_4(N))$  where
$N^{(2)}$ is the 2-skeleton of $N$.  
For any such choice of Spin structure,  $\alpha$ maps to $\kervaire^2(f) \cap [N]$ in
$E_{2,2}^\infty = H_2(N;\Omega^{Spin}_2)$.  

An analogous statement is true in the
topological category.
\end{theorem}
\begin{corollary}  \label{CSK} If $\kappa_2 : H_2(\pi;{\Bbb Z}_2) \ra L_4({\Bbb Z}\pi)$  is injective
and $M$ and $N$ are closed, homotopy equivalent, smooth (resp. topological) 4-manifolds 
which admit Spin structures and have fundamental groups isomorphic to $\pi$, then there
exist Spin structures $\tilde{\nu}_M$ and $\tilde{\nu}_N$, and maps $F_M : M \ra B\pi$,
$F_N : N \ra B\pi$ which induce isomorphisms on the fundamental group, so that
$[(M,  \tilde{\nu}_M),F_M] = [(N,  \tilde{\nu}_N),F_N] \in \Omega^{Spin}_4(B\pi)$
(resp. $\Omega^{TopSpin}_4(B\pi)).$  Hence $M$ and $N$ are stably diffeomorphic (resp. stably 
homeomorphic.)
\end{corollary}

\begin{proof}[Proof of Corollary] Let $f : M \to N$ be a homotopy equivalence.  First choose Spin structures $\tilde{\nu}_M$ and $\tilde{\nu}_N$ specified by the above theorem.  Then choose any map $F_N : N \to B\pi$ inducing an isomorphism on the fundamental group and let $F_M = F_N \circ f$.  Since we may assume the map $F_N$ is cellular, the above theorem shows
$$
\alpha = [(M,  \tilde{\nu}_M),F_M] - [(N,  \tilde{\nu}_N),F_N] \in F_{2,2}(\Omega^{Spin}_4(B\pi)).
$$
Since $\Omega_1^{Spin}(pt) = 0$ there is an split exact sequence
$$
0 \to \Omega_4^{Spin}(pt) \to F_{2,2} \to E^\infty_{2,2} = H_2(B\pi; \Omega^{Spin}_2)/\langle \text{im }d_2, \text{im }d_3 \rangle \to 0
$$

For a homotopy equivalence $f: M \to N$, the surgery obstruction $\theta(f)$ vanishes.  By 
the characteristic class formula (Proposition \ref{charclass}) and the injectivity of $\kappa_2$, one has $F_{N*}(
\kervaire^2(f) \cap [N]) = 0$.  Hence by the previous theorem, $\alpha$ maps to zero in $E^\infty_{2,2}$, and hence is detected by $\Omega^{Spin}_4$.  Since $M$ and $N$ have the same signature and since $\Omega^{Spin}_4$ and $\Omega^{TopSpin}_4$ are both detected by signature (see \cite[p. 325]{Kirby-Siebenmann}) we are done.
\end{proof}

\begin{proof}[Proof of Theorem]  Let $\hat{f} \in [N,G/O]$ be the classifying map for the degree one
normal map $f$.  We first lift $\hat{f}$ to $\hat{\alpha} \in [N,G]$, and then modify $\hat{\alpha}$
so that it becomes trivial on the 1-skeleton $N^{(1)}$.  Let $\xi$ be the bundle over $N$ given
as part of the degree one normal map structure.  Since $f$ is a degree one
normal map with $\sign N = \sign M$, it is easy to check that $w_2(\xi) = w_2(\nu_N)$
and $p_1(\xi) = p_1(\nu_N)$, and hence $\xi$ and $\nu_N$ are isomorphic bundles over the
4-manifold $N$.  Thus we may assume our map $\hat{f}$ classifies a fiber homotopy equivalence
$N \times S^k \ra N \times S^k$ which gives a map $\hat{\alpha}: N \ra G_{k+1} \subset G$.
By possibly multiplying by a matrix of determinant $-1$, we can assume that $\hat{\alpha}$
lands in the orientation preserving component $SG$.  Let $b \in H^1(N; \pi_1G)$ be the 
first obstruction to null-homotopy.  There is a map $\phi : \Bbb RP^\infty \ra
SO$ inducing an isomorphism on $\pi_1$.  Here $\phi([v]) = R_* \circ R_v$ where $R_w$ is
reflection through the line $\Bbb Rw$ and $* = (1,0,. . ., 0)$.  Thus by representing $b$ by a map to $\Bbb RP^\infty$, one can find
 an element $B \in [N,SO]$
whose first obstruction to null-homotopy is $b$ (using the identification $\pi_1O = \pi_1G)$.
Replacing the lift $\hat{\alpha}$ by $\hat{\alpha}-B$ we get our lift $\hat{\alpha} \in [N,G]$
of $\hat{f}$ so that $\hat{\alpha}(N^{(1)}) = *$.  Now $\hat{\alpha}$ lands in the 
orientation-preserving component $SG_k$ and the 
first obstruction to null-homotopy is $\kervaire^2(f)
\in H^2(N;\pi_2G)$.

To a map $g : N \ra SG_k$ assign the map $\Gamma(g):\Sigma^kN_+ \ra S^k$ given by
$$\Gamma(g) = ((Ad\ g) \vee (pr_1 \circ \nu_\Sigma))\circ\mu_\Sigma$$
where  $\mu_\Sigma$ is the comultiplication given by the first suspension coordinate
$\Sigma^kN_+$, $\nu_\Sigma$ is the coinverse, $pr_1$ is given by collapsing $N$ to a point
and $(Ad\ g)(t,n) = g(n)t$.  (In other words, $\Gamma(g) = Ad\ g - pr_1.$)
This gives a well-defined transformation
$\Gamma : [N,SG] \ra [N,\underline{S}^0]$ where  $\underline{S}^0$ is the sphere
spectrum.  If  $g(N^{(1)}) = *$, then using the coinverse,
$\Gamma(g)$ admits a canonical homotopy to a map trivial on
$\Sigma^kN^{(1)}_+$. We claim that the first obstruction for $g$ being null-homotopic
equals that of $\Gamma(g)$ in $H^2(N;\pi_2G)) = H^2(N;\pi_2\underline{S}^0)$.  This
follows from the identification of $[S^2,SG_k]$ with $[S^{k+2},S^k]$ which goes as
follows (see \cite{Milgram-Madsen}).  
\begin{eqnarray}
[S^2,SG_k]  & = & [(S^2_+,+),(SG_k,Id)]\nonumber\\
&\cong& [(S^2_+,+),(\Omega^k_1S^k,Id)]\nonumber\\ & \cong& [(S^2_+,+),(\Omega^k_0S^k,*)]\nonumber\\
&\cong& [(\Sigma^kS^2_+,+),(S^k,*)]\nonumber\\& = &[S^{k+2},S^k]\nonumber
\end{eqnarray}
where $\Omega^k_iS^k$ is the degree $i$ component of $\Omega^kS^k$, the first
isomorphism comes from the inclusion of $\Omega^k_1S^k \ra SG_k$ as basepoint
preserving degree 1 maps, the second isomorphism by $\mu_\Omega(- \wedge \nu_\Omega
\circ (Ad\ Id))$ and the third isomorphism is given by the adjoint.  Note that before the
adjoint correspondence can be applied, the basepoint needs to be the constant loop.

There is a map $i :\underline{S}^0 \ra \underline{MSpin}$ given by inclusion of a fiber,
which induces an isomorphism on  $\pi_i$ for $i \le 2$.  (See \cite{Milnor} for information 
on Spin structures.)  Thus the first obstruction to null-homotopy of $i\circ\Gamma(g) : 
\Sigma^kN_+ \ra MSpin(k)$ 
is still $\kervaire^2(f).$  On the other hand if $N$ is given a Spin structure $\tilde{\nu}_N$,
the Poincar\'e dual of $i\circ\Gamma(g)$ is given by $\alpha = [(M,\tilde{\nu}_M), f  ] -
[(N,\tilde{\nu}_N), Id  ]$. Hence by \ref{duality lemma} and 
\ref{obstruction lemma}, we have
identified images of $\alpha$ and  $\kervaire^2(f) \cap [N]$ in $E^\infty_{2,2}$.  

The uniqueness of the Spin structure on $N$
follows from the fact that a change in the Spin structure would make the Poincar\'e dual of $\alpha$
in $E^{1,-1}_\infty \subset H^1(N;\Omega^{-1}_{Spin}) = H^1(N; {\Bbb Z}_2)$ non-zero.

Finally, according to \cite[Lemma 5.3.2]{Teichner(92)}, $d_2 : E^{2}_{4,1}  \to E^2_{2,2}$ is the dual of $Sq^2 : H^2(N; \Z_2) \to H^4(N; \Z_2)$ which is zero since $N$ is Spin.  This shows that $E^\infty_{2,2} = H_2(N; \Omega_2^{Spin})$.  (This last fact was occurred in the statement of Theorem \ref{SK}, but was not used in the proof of  Corollary \ref{CSK}).
\end{proof}

For the non-Spin case on has:

\begin{theorem} \label{SKno} Let $f:M\ra N$  be a degree one normal map between closed, oriented, smooth 4-manifolds 
with the same signature inducing an isomorphism on fundamental group,
 and so that $w_2\widetilde{M}$ and $w_2\tilde{N}$ are both 
zero.  
Then there are normal
1-smoothings $\tilde{\nu}_M$ and $\tilde{\nu}_N$ in $\xi = \xi(N,0,w_2(\nu_N))$ 
 so that $\alpha = [M,\tilde{\nu}_M] -
[N, \tilde{\nu}_N]$ is in the filtration subgroup  $F_{2,2}
(\Omega_4(\xi)))$ of the James spectral sequence and $\alpha$ maps to $\kervaire^2(f) \cap [N]$ in
$E_{2,2}^\infty = H_2(N;(\Omega^{Spin}_2)^{w_1})$.  An analogous statement is true in the
topological category, provided that, in addition, $\ks (M) = \ks (N)$.
\end{theorem}

\begin{proof}  The proof runs parallel to the Spin case up to the point where we obtain
$i \circ \Gamma(g) \in [N_+,MSpin]$.  Next we use that $D(\nu_N)$ is canonical framed
and so admits a Spin structure.  We use Poincar\'e duality $\Omega^0_{Spin}(D(\nu_N)) \cong
\Omega^{Spin}_{k}(D(\nu_N),S(\nu_N))$ and see that $i \circ \Gamma(g)$ maps to
$[D(f^*(\nu_N))] - [D(\nu_N)]$, where both manifolds are framed and hence
have Spin structures.  To identify this with the element from the James 
spectral sequence we use Lemma \ref{James}.
\end{proof}

In the topological case, there is a lemma that may be of independent interest.

\begin{lemma}  If $\kappa_2$ is injective for $\pi$, and $f : M \ra N$
is a homotopy equivalence between closed, orientable, topological 4-manifolds 
whose fundamental groups are isomorphic to $\pi$ and whose universal covers admit
Spin structures, then $\ks (M) = \ks (N)$.
\end{lemma}

\begin{proof}
Since $f$ is a degree one normal map we have
$$\ks (M) - \ks (N) = (\sign(M) - \sign(N))/8 + \kervaire^2(f)^2$$
(see \cite[p. 328-330]{Kirby-Siebenmann}).
Since $f$ is a homotopy equivalence the signature term vanishes.  Let $\alpha = \kervaire^2(f)$,
let $w \in H^2(\pi;{\Bbb Z}_2)$ map to $w_2(N)$, and let $F : N \ra B\pi$ induce an
isomorphism on $\pi$.  Since $f$ is a homotopy equivalence
and $\kappa_2$ is injective, $F_*(\alpha \cap [N]) = 0$ by Proposition \ref{charclass}. We have
\begin{eqnarray}
\alpha^2  & = & (w_2 \cup \alpha)[N]\nonumber\\
& = & w_2(\alpha \cap [N])\nonumber\\ & = & w(F_*(\alpha \cap [N]))\nonumber\\
& = & 0\nonumber
\end{eqnarray}  
\end{proof}

\begin{corollary}
If $\kappa_2$ is injective for $\pi$, and $f : M \ra N$
is a homotopy equivalence between orientable, topological 4-manifolds whose universal cover is Spin
and with fundamental group $\pi$,
then $M$ and $N$ are stably homeomorphic.
\end{corollary}

As advertised.

\begin{example} \label{FKS}
 Let $N$ be a closed, orientable, smooth 4-manifold with fundamental group 
$\Bbb Z_2$, non-trivial $w_2$, but so that $w_2(\tilde{N}) = 0$.  (The existence of 
such an $N$ is shown by
computing $\Omega_4(\xi(B\Bbb Z_2,0,\not= 0))$ or is explicitly given by
a quotient of a free involution on the Kummer surface.)  Then there is an
$\alpha \in H^2(N; \Bbb Z_2)$ so that $\alpha^2 \not= 0$.  Choose $\hat{f}
\in [N,G/TOP] = H^4(N; \Bbb Z) \oplus H^2(N; \Bbb Z_2)$ corresponding to
$(0,\alpha)$.  Since $\kappa_2$ vanishes for $\pi = \Bbb Z_2$ (see 
\cite[Prop. 7.4]{Hambleton-Milgram-Taylor-Williams}), there is no obstruction to surgery to
a homotopy equivalence $M \ra N$. Then $\ks(M) = \alpha^2 \not= 0$.
Thus the Kirby-Siebenmann invariant is not homotopy invariant for closed 4-manifolds whose universal cover is Spin.  (Note the fundamental group in this example is not torsion-free.) 
\end{example}

\noindent {\bf Note:}  The paper of Teichner \cite{Teichner(97)}, dealing with homotopy invariance of the Kirby-Siebenmann invariant for 4-manifolds, could be profitably revisited from the point of view of the homomorphism $I_0 \oplus \kappa_2$. 

The survey of Kirby and Taylor \cite[Theorem 22]{Kirby-Taylor} gives quite a different approach to Theorem \ref{main} of our paper.  It may be interesting to compare the two approaches.

\end{document}